\documentclass[12pt]{article}
\usepackage{amsfonts,amsmath,amsxtra}
\usepackage{amssymb,amscd}
\usepackage[mathscr]{eucal}
\def\hybrid{\topmargin 0pt      \oddsidemargin 0pt
        \headheight 0pt \headsep 0pt
        \textwidth 16.5cm
        \textheight 23cm
        \marginparwidth 0.0in
        \parskip 5pt plus 1pt   \jot = 1.5ex}
\catcode`\@=11
\def\marginnote#1{}
\newcount\hour
\newcount\minute
\newtoks\amorpm
\hour=\time\divide\hour by60
\minute=\time{\multiply\hour by60 \global\advance\minute by-\hour}
\edef\standardtime{{\ifnum\hour<12 \global\amorpm={am}%
        \else\global\amorpm={pm}\advance\hour by-12 \fi
        \ifnum\hour=0 \hour=12 \fi
      \number\hour:\ifnum\minute<10 0\fi\number\minute\the\amorpm}}
\edef\militarytime{\number\hour:\ifnum\minute<10 0\fi\number\minute}

\def\draftlabel#1{{\@bsphack\if@filesw {\let\thepage\relax
   \xdef\@gtempa{\write\@auxout{\string
      \newlabel{#1}{{\@currentlabel}{\thepage}}}}}\@gtempa
   \if@nobreak \ifvmode\nobreak\fi\fi\fi\@esphack}
        \gdef\@eqnlabel{#1}}
\def\@eqnlabel{}
\def\@vacuum{}
\def\draftmarginnote#1{\marginpar{\raggedright\scriptsize\tt#1}}

\def\draft{\oddsidemargin -0.1truein
        \def\@oddfoot{\sl preliminary draft \hfil
        \rm\thepage\hfil\sl\today\quad\militarytime}
        \let\@evenfoot\@oddfoot \overfullrule 3pt
        \let\label=\draftlabel
        \let\marginnote=\draftmarginnote
\def\@eqnnum{{\rm (\theequation)}
\rlap{\kern\marginparsep\tt\@eqnlabel}%
\global\let\@eqnlabel\@vacuum}  }

\newfont{\Bbbb}{msbm7 scaled 1\@ptsize00}
\newcommand{\zs}{\raise-1pt\hbox{$\mbox{\Bbbb Z}$}}

\font\sevenmsa=msam6 
\newfam\msafam
\textfont\msafam=\sevenmsa
\def\hexnumber@#1{\ifnum#1<10 \number#1\else
\ifnum#1=10 A\else\ifnum#1=11 B\else\ifnum#1=12 C\else
\ifnum#1=13 D\else\ifnum#1=14 E\else\ifnum#1=15 F\fi\fi\fi\fi\fi\fi\fi}
\def\msa@{\hexnumber@\msafam}
\def\llcorner{\delimiter"4\msa@78\msa@78 }
\def\lrcorner{\delimiter"5\msa@79\msa@79 }
\mathchardef\blacktriangleright="3\msa@49
\mathchardef\blacktriangleleft="3\msa@4A
\font\tenmsb=msbm10 scaled 1\@ptsize00
\newfam\msbfam
\textfont\msbfam=\tenmsb
\scriptfont\msbfam=\tenmsb

\newdimen\linethick  \linethick=0.4pt
\newdimen\hboxitspace    \hboxitspace=5pt
\newdimen\vboxitspace    \vboxitspace=5pt

\def\fr#1{%
\be\new
\vcenter{
\hrule height\linethick
           \hbox{\vrule width\linethick
                 \kern\hboxitspace
                 \vbox{\kern\vboxitspace
                       \hbox{$\begin{array}{c}\displaystyle#1
          \end{array}$}%
                       \kern\vboxitspace}%
                 \kern\hboxitspace
                 \vrule width\linethick}%
           \hrule height\linethick}%
\ee}

\newdimen\Squaresize \Squaresize=14pt
\newdimen\Thickness \Thickness=0.5pt

\def\Square#1{\hbox{\vrule width \Thickness
   \vbox to \Squaresize{\hrule height \Thickness\vss
      \hbox to \Squaresize{\hss#1\hss}
   \vss\hrule height\Thickness}
\unskip\vrule width \Thickness}
\kern-\Thickness}

\def\Vsquare#1{\vbox{\Square{$#1$}}\kern-\Thickness}

\def\numberbysection{\@addtoreset{equation}{section}
        \def\theequation{\thesection.\arabic{equation}}}
\numberbysection

\renewcommand{\theequation}{\thesection.\arabic{equation}}
\def\titlepage{\@restonecolfalse\if@twocolumn\@restonecoltrue\onecolumn
     \else \newpage \fi \thispagestyle{empty}\c@page\z@
        \def\thefootnote{\fnsymbol{footnote}} }

\def\endtitlepage{\if@restonecol\twocolumn \else  \fi
        \def\thefootnote{\arabic{footnote}}
        \setcounter{footnote}{0}}  
\relax

\hybrid
\makeatletter
\newdimen\normalarrayskip            
\newdimen\minarrayskip               
\normalarrayskip\baselineskip
\minarrayskip\jot
\newif\ifold             \oldtrue            \def\new{\oldfalse}
\def\arraymode{\ifold\relax\else\displaystyle\fi}
\def\eqnumphantom{\phantom{(\theequation)}} 
\def\@arrayskip{\ifold\baselineskip\z@\lineskip\z@
     \else
     \baselineskip\minarrayskip\lineskip1\baselineskip\fi}

\def\@arrayclassz{\ifcase \@lastchclass \@acolampacol \or
\@ampacol \or \or \or \@addamp \or
   \@acolampacol \or \@firstampfalse \@acol \fi
\edef\@preamble{\@preamble
  \ifcase \@chnum
     \hfil$\relax\arraymode\@sharp$\hfil
     \or $\relax\arraymode\@sharp$\hfil
     \or \hfil$\relax\arraymode\@sharp$\fi}}

\def\@array[#1]#2{\setbox\@arstrutbox=\hbox{\vrule
     height\arraystretch \ht\strutbox
     depth\arraystretch \dp\strutbox
width\z@}\@mkpream{#2}\edef\@preamble{\halign \noexpand\@halignto
\bgroup \tabskip\z@ \@arstrut \@preamble \tabskip\z@ \cr}%
\let\@startpbox\@@startpbox \let\@endpbox\@@endpbox
  \if #1t\vtop \else \if#1b\vbox \else \vcenter \fi\fi
  \bgroup \let\par\relax
  \let\@sharp##\let\protect\relax
  \@arrayskip\@preamble}
\def\eqnarray{\stepcounter{equation}%
              \let\@currentlabel=\theequation
              \global\@eqnswtrue
              \global\@eqcnt\z@
              \tabskip\@centering              
              \let\\=\@eqncr
              $$%
            \halign to \displaywidth  \bgroup
             \eqnumphantom \@eqnsel
      \hskip\@centering                               
    $\displaystyle  \tabskip\z@ {##}$%
    &\global\@eqcnt\@ne \hskip 2\arraycolsep
         $ \displaystyle  \arraymode{##}$\hfil
    &\global\@eqcnt\tw@ \hskip 2\arraycolsep
         $\displaystyle\tabskip\z@{##}$\hfil
         \tabskip\@centering
    &{##}\tabskip\z@\cr}
\makeatother

\newcommand{\beq}[1]{\begin{equation}\label{#1}}
\newcommand\eeq{\end{equation}}
\newcommand\bqa{\begin{eqnarray}}
\newcommand\eqa{\end{eqnarray}}
\def\be{\begin{eqnarray}\new\begin{array}{cc}}
\def\ee{\end{array}\end{eqnarray}}

\def\beq{\begin{equation}}
\def\eeq{\end{equation}}
\def\bse{\begin{subequations}}                
\def\ese{\end{subequations}}
\def\bp{\begin{pmatrix}}
\def\ep{\end{pmatrix}}

\def\stack#1#2{\raise0.7pt\hbox{$\mathrel{\mathop{#2}\limits^{#1}}$}}
\def\tr{\triangleright}
\def\tl{\triangleleft}
\def\sem{\mathsurround=0pt \raise1pt
\hbox{$\scriptscriptstyle>\!\!$}\:\!\!\tl}
\def\mes{\mathsurround=0pt \tr\!\:\!\raise0.8pt
\hbox{$\scriptscriptstyle\!\!<$}\,}
\def\]{\mathsurround=0pt ]\raise-2pt\hbox{$_\ast$}}

\def\<{\langle}
\def\>{\rangle}

\def\CH{\mathcal{H}}

\def\we{\raise-1pt\hbox{$\,\stackrel{\wedge}{,}\,$}}

0
\begin{document}

\footnotesize
\normalsize

\newpage

\thispagestyle{empty}

\begin{center}

\phantom.
\bigskip
{\hfill{\normalsize hep-th/0608152}\\
\hfill{\normalsize ITEP-TH-XX/06}\\
\hfill{\normalsize HMI-06-XX}\\
\hfill{\normalsize TCD-MATH-06-XX}\\
[10mm]\Large\bf
Givental Integral Representation  for Classical Groups
\footnote{To be published in the Proceedings of the Satellite ICM 2006
conference: "Integrable systems in Applied Mathematics",
Colmenarejo (Madrid, Spain), 7-12 September 2006.}}
\vspace{0.5cm}

\bigskip\bigskip
{\large A. Gerasimov}
\\ \bigskip
{\it Institute for Theoretical and
Experimental Physics, 117259, Moscow,  Russia} \\
{\it  School of Mathematics, Trinity
College, Dublin 2, Ireland } \\
{\it  Hamilton
Mathematics Institute, TCD, Dublin 2, Ireland},\\
\bigskip
{\large D. Lebedev\footnote{E-mail: lebedev@ihes.fr}}
\\ \bigskip
{\it Institute for Theoretical and Experimental Physics, 117259,
Moscow, Russia}  \\
{\it l'Institute des Hautes
\'Etudes Scientifiques, 35 route de Chartres,
Bures-sur-Yvette, France},\\
\bigskip
{\large S. Oblezin} \footnote{E-mail: Sergey.Oblezin@itep.ru}\\
\bigskip {\it Institute for Theoretical and Experimental Physics,
117259, Moscow,
Russia}\\
{\it Max-Planck-Institut f\"ur Mathematik, Vivatsgasse 7, D-53111
Bonn, Germany},\\
\end{center}

\vspace{0.5cm}

\begin{abstract}
\noindent

We propose  integral representations for  wave functions of 
$B_n$, $C_n$, and $D_n$ open Toda chains  at zero eigenvalues
of the Hamiltonian operators   
thus  generalizing  Givental representation for $A_n$.  We also
construct Baxter $Q$-operators for closed Toda chains corresponding
to Lie algebras $B_{\infty}$, $C_{\infty}$, $D_{\infty}$, affine Lie
algebras  $B^{(1)}_n$, $C^{(1)}_n$, $D^{(1)}_n$ and twisted affine Lie
algebras $A^{(2)}_{2n-1}$ and $A^{(2)}_{2n}$. Our approach is
based on a generalization of  the  connection between  Baxter $Q$-operator for
$A_n^{(1)}$ closed Toda chain and Givental representation
for the wave function of $A_n$ open Toda chain uncovered previously.

\end{abstract}

\vspace{1cm}

\clearpage \newpage


\normalsize
\section{Introduction}

A remarkable  integral representation  for  the common
eigenfunctions of $A_n$ open Toda chain  Hamiltonian operators
was proposed in  \cite{Gi} (see also
\cite{JK}). This representation is based on  
a flat  degeneration of $A_n$ flag manifolds
to a Gorenstein toric Fano variety   
(see \cite{L},\cite{Ba}, \cite{BCFKS} for details).  
This results  in a purely combinatorial 
description of the integrand in  the integral representation
  \cite{Gi}. An important application of the Givental  
integral representation so far was an explicit construction 
of the mirror dual of $A_n$ flag manifolds. 

Later it turns out  that the representation introduced in \cite{Gi}
is also interesting from another points of view. Thus  
it was shown in  \cite{GKLO}  that this
 integral representation has  natural iterative structure allowing
the connection of  $A_{n-1}$ and $A_{n}$ wave functions by a  simple integral
transformation. It was demonstrated  that thus defined integral
transformation is given by a degenerate version of the Baxter
$Q$-operator realizing quantum B\"acklund transformations in
 closed Toda chain \cite{PG}. Let us note that the torification 
of  flag manifolds leads to a distinguished set of coordinates 
 on its open parts. A group theory construction of these coordinates 
is also connected with a degenerate $Q$-operator and was clarified  
in \cite{GKLO}.

Up to now  the Givental integral representation  
was only generalized \cite{BCFKS} to the case of  
degenerate $A_n$ open Toda chains \cite{STS}  leading to a construction of the 
mirror duals  to  partial flag manifolds $G/P$ for $G=SL(n+1,\mathbb{C})$ and 
$P$ being a parabolic subgroup. A natural approach 
to  a generalization of the integral representation 
to other Lie algebras  could be based 
on  the relation with  Baxter 
$Q$-operator. However no  generalization of the
Baxter $Q$-operator to other Lie algebras  was  known. 
In this note  we solve both these problems simultaneously
 for all  classical series of Lie algebras.   
We propose  a generalization of the Givental integral representation 
to  other classical series $B_n$, $C_n$, $D_n$ 
and construct $Q$-operators for  affine Lie algebras
$A_n^{(1)}$, $A^{(2)}_{2n}$, $A_{2n-1}^{(2)}$, 
$B_n^{(1)}$, $C_n^{(1)}$, $D_n^{(1)}$ and infinite Lie algebras 
$B_{\infty}$, $C_{\infty}$ and $D_{\infty}$.  
We also generalize the connection between $Q$-operators and 
integral representations of wave functions 
 to  all classical series. 

Let us stress that 
there is an important difference in the construction of the integral representations
between $A_n$   and  other classical series. 
The  kernel of the integral operator providing
recursive  construction of the integral representation for $A_n$  has
 a simple form of the exponent of the sum of  exponents in the natural
coordinates. For other  classical groups
recursive  operators of the same type   exist  but they 
relate   Toda chain  wave functions for  different classical series (e.g.
 $C_{n}$ and $D_{n}$). Integral
operators connecting Toda wave functions  in
the same series (e.g. $D_n$ and $D_{n-1}$) are given by  compositions of the
elementary integral operators. 

In this note we restrict ourselves by   explicit constructions of the integral
representations of eigenfunctions  of the quadratic open Toda chain 
Hamiltonian operators at zero eigenvalues. The general case 
of all Hamiltonians and non-zero eigenvalues 
will be published elsewhere. Also we leave for another 
occasion the elucidation of  a group theory interpretation of the obtained 
results. 

The plan of this paper  is as follows. In Section 2 we summarize
the results of \cite{GKLO}. In Section 3
we construct  kernels of the elementary
integral operators intertwining  Hamiltonian
operators of open Toda chains for (in general different)  
classical series of finite Lie algebras.
In Section 4 using  the results from Section 3 we give explicit
integral representations for the wave functions of
$B_n$, $C_n$ and $D_n$ open Toda chains.
In Section 5 we describe a generalization of Givental  
diagrams to other classical series and remark 
on the connection with toric degeneration 
of $B_n$, $C_n$ and $D_n$  flag manifolds. 
In Section 6 we construct elementary integral operators intertwining
 Hamiltonian operators  of  closed Toda chains for  (in general
 different) classical series of  affine Lie algebras.
In Section 7   we construct  integral kernels
for Baxter $Q$-operators for all classical series of (twisted)
 affine Lie algebras as appropriate compositions of the elementary
intertwining operators. In Section 8 we construct 
$Q$ operators for $B_{\infty}$, $C_{\infty}$ and $D_{\infty}$ Toda  
chains.  We conclude in Section 9 with a short
discussion of the results presented in this note. 

We were informed by E.~Sklyanin
that he also has some progress in the construction of Baxter
$Q$-operators for Toda theories.

{\em Acknowledgments}:  The authors are grateful to S.~Kharchev for
 discussions at the initial stage of this project  and  to 
B.~Dubrovin and M.~Kontsevich
for their interest in this work.  The research of A.~Gerasimov
was  partly supported by the Enterprise Ireland Basic
Research Grant. D.~Lebedev is grateful to
 Institute des Hautes \'Etudes Scientifiques for  warm
hospitality. S.~Oblezin is grateful to Max-Planck-Institut f\"ur Mathematik
for excellent working conditions.

\section{Recursive structure of Givental representation}

In this section we recall a recursive construction of the Givental
integral representation discussed in \cite{GKLO}.

The solution of a quantum integrable system starts with the
finding of the full set of common eigenfunctions of the quantum
Hamiltonian operators of $A_n$ Toda chain. 
 Note that the difference between wave functions 
for $\mathfrak{sl}_n$ and $\mathfrak{gl}_n$ manifests only at non-zero
eigenvalues of the  element of the center of $U\mathfrak{gl}_n$
linear over the generators. 
In the following we will consider 
only the wave functions corresponding to zero 
eigenvalues of all elements of the center. 
Thus we will  always consider $\mathfrak{gl}_n$ Toda chain
instead of $A_n$ Toda chain. 
The quadratic quantum Hamiltonian of 
$\mathfrak{gl}_n$ open Toda chain is given by \be
H^{\mathfrak{gl}_n}(x)=-\frac{\hbar^2}{2}\sum\limits_{i=1}^{n}\frac{\partial^2}{\partial
x_i^2}+ \sum\limits_{i=1}^{n-1}g_i e^{x_{i+1}-x_{i}}\,. \ee
In \cite{Gi} the following remarkable representation
 for a  common eigenfunction of quantum Hamiltonians
of $\mathfrak{gl}_n$ open Toda chain was proposed
\be\label{intrep} \Psi(x_1,\ldots,x_{n})=\int_{\Gamma}
e^{\frac{1}{\hbar}\mathcal{F}_n(x)}\bigwedge_{k=1}^{n-1} \bigwedge_{i=1}^k
dx_{k,i}, \ee where $x_{n,i}:=x_i$, the function $\mathcal{F}_n(x)$
is given by \be\label{pot} \hspace{-0.5cm}
\mathcal{F}_n(x)=\sum_{k=1}^{n-1}\sum_{i=1}^k\Big(e^{x_{k,i}-x_{k+1,i}}
+g_ie^{x_{k+1,i+1}-x_{k,i}}\Big), \ee and the cycle $\Gamma$ is a
middle dimensional submanifold in the $n(n-1)/2$- dimensional
complex torus with coordinates $\{\exp\,{x_{k,i}},\,
i=1,\ldots,k;\,k=1,\ldots, n-1\}$ such that the integral converges.
The eigenfunction (\ref{intrep}) solves the  equation 
\bqa H^{\mathfrak{gl}_n}(x)\,\Psi(x_1,\cdots, x_n)= 0.  \eqa In the
following we put $\hbar=1$ for convenience. 

The derivation of the integral
representation (\ref{intrep}) using the recursion over the rank  $n$
of the Lie algebra $\mathfrak{gl}_n$  was given in \cite{GKLO}.
The integral representation for the wave function 
can be represented in the following form
\bqa\label{recintrep} \Psi(x_{1},\ldots,x_{n})\,=\,
\int\bigwedge_{k=1}^{n-1}\bigwedge_{i=1}^kdx_{k,i}\,
\prod_{k=1}^{n-1}Q_{k+1,\,k}(x_{k+1,1},\ldots,x_{k+1,k+1};
\,x_{k,1},\ldots,x_{k,k}), \eqa with the integral kernel
\be\label{QBAXTER} Q_{k+1,k}(x_{k+1,i};\, x_{k,i})=
\exp\Big\{\, \sum_{i=1}^{k}
\left(e^{x_{k,i}-x_{k+1,i}}+g_ie^{x_{k+1,i+1}-x_{k,i}}\right)\,\Big\}. \ee Here
we have $x_i:=x_{n,i}$. The following differential equation 
for the kernel holds  \bqa\label{intertw}
H^{\mathfrak{gl}_{k+1}}(x_{k+1,i})Q_{k+1,k}(x_{k+1,i},\, x_{k,i})=
Q_{k+1,k}(x_{k+1,i},\, x_{k,i})\,H^{\mathfrak{gl}_{k}}(x_{k,i}), \eqa where \be
H^{\mathfrak{gl}_k}(x_i)
=-\frac{1}{2}\sum\limits_{i=1}^{k}\frac{\partial^2}{\partial
x_{k,i}^2}+ \sum\limits_{i=1}^{k-1}g_i e^{x_{k,i+1}-x_{k,i}}\, . \ee
Here and in the following we assume that in the relations similar to
(\ref{intertw}) the Hamiltonian operator on l.h.s. acts on the right and
the Hamiltonian on r.h.s. acts on the left. Thus the integral
operator with the kernel  $Q_{k+1,k}$ intertwines
 Hamiltonian operators for $\mathfrak{gl}_{k+1}$ and $\mathfrak{gl}_k$ open Toda chains.

The integral operator defined by the kernel (\ref{QBAXTER}) is
closely related with a Baxter $Q$-operator realizing B\"{a}cklund
transformations in a  closed Toda chain corresponding to affine Lie
algebra $\widehat{\mathfrak{gl}}_n$. Baxter $Q$-operator for zero
spectral parameter can be
written in the integral form with the   kernel \cite{PG}
\bqa\label{baff} Q^{\widehat{\mathfrak{gl}}_n}(x_i,y_i)=\exp\Big\{\,
\sum_{i=1}^{n} \left(e^{x_{i}-y_{i}}+g_ie^{y_{i+1}-x_{i}}\right)\,\Big\}, \eqa
where $x_{i+n}=x_i$ and $y_{i+n}=y_i$. This operator commutes with
the Hamiltonian operators of the closed Toda chain. Thus for example
for the quadratic Hamiltonian we have \bqa\label{intertwaff}
\CH^{\widehat{\mathfrak{gl}}_n}(x_i)Q^{\widehat{\mathfrak{gl}}_n}(x_i,\,
y_i)= Q^{\widehat{\mathfrak{gl}}_n} (x_i,\,
y_i)\CH^{\widehat{\mathfrak{gl}}_n}(y_i), \eqa where  \be
\CH^{\widehat{\mathfrak{gl}}_n}=
-\frac{1}{2}\sum\limits_{i=1}^{n}\frac{\partial^2}{\partial x_i^2}+
\sum\limits_{i=1}^{n}g_i e^{x_{i+1}-x_{i}}\,. \ee Here we impose the 
conditions  $x_{i+n}=x_i$. The recursive operator (\ref{QBAXTER}) can be
obtained from Baxter operator (\ref{baff}) in the limit
$g_{n}\rightarrow 0$, $x_n\rightarrow-\infty$.

The main objective of this note is to generalize  the
representation  (\ref{intrep}), (\ref{pot}) to  other
classical series   $B_n$, $C_n$ and $D_n$
of finite Lie groups. Before we present the 
integral representations  for $B_n$, $C_n$ and $D_n$
 let us comment on the main subtlety in their constructions.  
As it was
explained in \cite{GKLO} the variables $x_{k,i}$ in the integral
representation for $A_n$ have a clear meaning of the linear coordinates on Cartan
subalgebras of the intermediate Lie algebras entering recursion
 $A_n\to A_{n-1}\to \cdots \to A_{1}$.   This is a consequence of the
identity
\be ({\rm dim}\, (\mathfrak{gl}_n)-{\rm rk}\,(\mathfrak{gl}_n)) -({\rm
  dim}\, (\mathfrak{gl}_{n-1})-{\rm rk}\,(\mathfrak{gl}_{n-1}))
=2\,{\rm rk}\, \mathfrak{gl}_{n-1}. \ee
However for  other classical series
there is no such simple relation. In general one finds 
more integration variables in the integral representation then
those  arising as  linear coordinates on  intermediate Cartan
subalgebras. It turns out that the 
elementary integral operators  intertwine Hamiltonians corresponding
to  Toda chains for  {\it different} Lie algebras. The recursive operators
are then constructed as appropriate 
 compositions of the elementary intertwining operators.

\section{Elementary intertwiners for open Toda chains}

Let $\mathfrak{g}$ be
a simple Lie algebra, $\mathfrak{h}$ be a  Cartan
subalgebra, $n=\dim \mathfrak{h}$ be the rank of $\mathfrak{g}$,
$R\subset \mathfrak{h}^*$  be the root system,  $W$ be  the  Weyl
group. Let us fix a decomposition $R=R_+\cup R_-$ of the roots on
positive and negative roots. Let
$\alpha_1,\ldots,\alpha_n$ be the bases of  simple roots.
Let $(,)$ be a $W$-invariant bilinear symmetric form on
$\mathfrak{h}^*$ normalized so that $(\alpha,\alpha)=2$ for a long root.
This form provides an identification of $\mathfrak{h}$
with $\mathfrak{h}^*$ and thus can be considered as a bilinear form on
$\mathfrak{h}$. Choose an orthonormal basis $e=\{e_1,\ldots,e_{n}\}$
in $\mathfrak{h}$. Then for any $x\in \mathfrak{h}$
one has  a decomposition $x=\sum_{i=1}^n x_ie_i$.
One associates  with these data an open
Toda chain with a quadratic Hamiltonian \bqa H^{R}(x_i)=
-\frac{1}{2}\sum_{i=1}^{n}\frac{\partial^2}{\partial x_i^2}\,+\,
\sum_{i=1}^ng_ie^{\alpha_i(x)}\, .\eqa   For the standard facts on  Toda theories
corresponding to arbitrary root systems see e.g. \cite{RSTS}.

We start with explicit expressions for  elementary intertwining 
operators for  open Toda chains. The necessary facts on the root
systems (including non-reduced ones) can be found in \cite{He}.

\subsection{ $BC\leftrightarrow  B$}

Let $e=\{e_1,\ldots,e_{n}\}$ be an  orthonormal basis in
$\mathbb{R}^n$. Non-reduced root system of type $BC_n$ can be
defined as \bqa \alpha_0=2e_1,\qquad \alpha_1=e_1,\qquad
\alpha_{i+1}=e_{i+1}-e_i,\qquad 1\leq i\leq n-1 \eqa and the
corresponding Dynkin diagram is \bqa
\frac{\alpha_0}{\alpha_1}\,\Longleftrightarrow\,\alpha_2
\longleftarrow\ldots\longleftarrow\alpha_n \eqa where the first
vertex from the left is a doubled vertex corresponding to a   
 reduced $\alpha_1=e_1$ and non-reduced
 $\alpha_0=2e_1=2\alpha_1$ roots.

Quadratic Hamiltonian operator of the corresponding open Toda chain
is given by \bqa H^{BC_n}(x_i)=-
\frac{1}{2}\sum_{i=1}^n\frac{\partial^2}{\partial x_i^2}+
\frac{g_1}{2}\Big(e^{x_1}+g_1e^{2x_1}\Big)+
\sum_{i=1}^{n-1}g_{i+1}e^{x_{i+1}-x_i}. \eqa Let us stress that the same
open Toda chain can be considered as a most general form of $C_n$
 open Toda chain  (see e.g. \cite{RSTS}, Remark p.61). However
in the following we will use the term $BC_n$ open Toda chain
to distinguish it from a more standard $C_n$ open Toda chain that
will be  consider below.

The  root system of type $B_n$ can be defined as \bqa \alpha_1=e_1,\qquad
\alpha_{i+1}=e_{i+1}-e_i,\qquad 1\leq i\leq n-1 \eqa and the
corresponding Dynkin diagram is \bqa
\alpha_1\,\Longleftarrow\,\alpha_2
\longleftarrow\ldots\longleftarrow\alpha_n\,. \eqa
Quadratic Hamiltonian operator of the corresponding open Toda chain
is given by \bqa H^{B_n}(x_i)=-
\frac{1}{2}\sum_{i=1}^{n}\frac{\partial^2}{\partial x_i^2}+
g_1e^{x_1}+\sum_{i=1}^{n-1}g_{i+1}e^{x_{i+1}-x_i}.\eqa An elementary
operator intertwining open Toda chain Hamiltonians for
 $BC_n$ and $B_{n-1}$ can be written in the integral form with
 the  kernel \be
Q_{BC_n}^{\,\,\,\,B_{n-1}}(z_1,\ldots,z_n;\,x_1,\ldots,x_{n-1})=
\exp\Big\{\,g_1e^{z_1}+\sum_{i=1}^{n-1}(e^{x_i-z_i}+
g_{i+1}e^{z_{i+1}-x_i})\,\Big\},\ee satisfying the following relation
 \bqa\label{BCB} H^{BC_n}(z)\,Q_{BC_n}^{\,\,\,\,B_{n-1}}(z,\, x)=
Q_{BC_n}^{\,\,\,\,B_{n-1}}(z,\, x) \,H^{B_{n-1}}(x) .
\eqa Similarly an elementary operator
intertwining $B_n$ and $BC_{n}$
 Hamiltonians  has an integral kernel
\be Q_{B_n}^{\,\,\,\,BC_{n}}(x_1,\ldots,x_n;\,z_1,\ldots,z_{n})=\\
\exp\Big\{\,g_1e^{z_1}+\sum_{i=1}^{n-1}\Big(e^{x_i-z_i}+
g_{i+1}e^{z_{i+1}-x_i}\Big)+e^{x_n-z_n}\,\Big\}.\ee

\subsection{ $C \leftrightarrow  D$}

The  root system of type $C_n$ can be defined as \bqa
 \alpha_{i}=e_{i+1}-e_i,\qquad \alpha_n=2e_n,\qquad
1\leq i\leq n-1,
\eqa and the   corresponding Dynkin diagram is \bqa
\alpha_1\,\longleftarrow\,\ldots \longleftarrow \alpha_{n-1}
\Longleftarrow\alpha_n\,. \eqa
Quadratic Hamiltonian operator of the corresponding open Toda chain
is given by
\be
H^{C_n}(x_i)=-
\frac{1}{2}\sum_{i=1}^{n}\frac{\partial^2}{\partial
x_i^2}+\sum_{i=1}^{n-1}g_{i}e^{x_{i+1}-x_i}+ 2g_ne^{-2x_n}.
\ee
The   root system of type $D_n$ is
\bqa \alpha_{i}=e_{i+1}-e_i,\qquad \alpha_n=-e_{n-1}-e_{n},\qquad
1\leq i< n, \eqa and the corresponding Dynkin diagram is
\bqa\begin{CD}
\alpha_1 @>>> \ldots @>>> \alpha_{n-2} @>>> \alpha_{n-1}\\
@.@. @VVV  @.\\ @.@. \alpha_{n}  @.\end{CD}\eqa
Quadratic Hamiltonian operator of the $D_n$ open Toda chain is given by\be
 H^{D_n}(x_i)=-
\frac{1}{2}\sum_{i=1}^n\frac{\partial^2}{\partial x_i^2}+
\sum_{i=1}^{n-1}g_ie^{x_{i+1}-x_i}+g_{n-1}g_ne^{-x_n-x_{n-1}}. \ee
An integral operator intertwining $C_n$ and $D_{n}$ Hamiltonians has
the kernel \be
Q_{C_n}^{\,\,\,\,D_{n}}(z_1,\ldots,z_n;x_1,\ldots,x_{n})=\\
\exp\Big\{\,
\sum_{i=1}^{n-1}\Big(e^{x_i-z_i}+g_{i}e^{z_{i+1}-x_i}\Big)+
e^{x_n-z_{n}}+g_ne^{-x_n-z_{n}}\Big\}.\ee Similarly an integral
operator with the kernel \be Q_{D_n}^{\,\,\,C_{n-1}}(x_1,\ldots,
x_n;\, z_1,\ldots,z_{n-1})=\\=
\exp\Big\{\sum_{i=1}^{n-1}\Big(e^{z_i-x_i}+g_ie^{x_{i+1}-z_i}\Big)+
g_ne^{-x_n-z_{n-1}}\Big\},\ee intertwines the following $D_n$ and
$C_{n-1}$ quadratic Hamiltonians \bqa H^{D_n}(x_i)=-
\frac{1}{2}\sum_{i=1}^n\frac{\partial^2}{\partial x_i^2}+
\sum_{i=1}^{n-1}g_ie^{x_{i+1}-x_i}+g_ne^{-x_n-x_{n-1}},\\
H^{C_{n-1}}(z_i)=-\frac{1}{2}\sum_{i=1}^{n-1}\frac{\partial^2}{\partial
z_i^2}+\sum_{i=1}^{n-2}g_{i}e^{z_{i+1}-z_i}+
2g_{n-1}g_ne^{-2z_{n-1}}.\eqa

\section{Givental representation for wave functions}

In the previous section  explicit expressions for  the kernels
 of elementary intertwining operators were presented. 
Now integral representations for  eigenfunctions of  open Toda chain
Hamiltonians are  given by a
quite straightforward generalization of $A_n$ case. Below we provide
integral representations for all classical series. For simplicity 
we put $g_i=1$ below. 

\subsection{$B_n$}
The eigenfunction for $B_n$ open Toda chain is given by
\be
\Psi^{B_n}(x_1,\ldots,x_n)\,=\,
\int\bigwedge_{k=1}^{n-1}\bigwedge_{i=1}^kdx_{k,i}\,
\prod_{k=1}^{n-1}Q_{B_{k+1}}^{\,\,B_k}(x_{k+1,1},\ldots,x_{k+1,k+1};
\,x_{k,1},\ldots,x_{k,k}),
\ee
where $x_i:=x_{n,i}$ and the kernels $Q_{B_{k+1}}^{\,\,\,\,B_k}$
of the integral operators are   given by the
convolutions of the kernels $Q_{B_{k+1}}^{\,\,\,\,BC_{k+1}}$
and $Q_{BC_{k+1}}^{\,\,\,\,B_k}$
\bqa
Q_{B_{k+1}}^{\,\,\,\,B_k}(x_{k+1,i};\,x_{k,j})=\int\bigwedge_{i=1}^kdz_{k,i}
\, Q_{B_{k+1}}^{\,\,\,\,BC_{k+1}}(x_{k+1,1},\ldots,x_{k+1,k+1};
z_{k+1,1},\ldots,z_{k+1,k+1})\times \\ \nonumber
\times  Q_{BC_{k+1}}^{\,\,\,\,B_k}
(z_{k+1,1},\ldots,z_{k+1,k+1};x_{k,1},\ldots,x_{k,k}).\eqa
Notice that the wave function is given by the integral over a contour 
of the real dimension equal to a complex  dimension of the flag manifold
$X=G/B$, where $G=SO(2n+1,\mathbb{C})$ and $B$ is a Borel subgroup   
 \be
\sum_{k=1}^n (2k-1)=n^2=|R_+|.
\ee

\subsection{ $C_n$}

The eigenfunction for $C_n$ open Toda chain is given by
\be
\Psi^{C_n}(z_1,\ldots,z_n)\,=\,
\int\bigwedge_{k=1}^{n-1}\bigwedge_{i=1}^kdz_{k,i}\,
\prod_{k=1}^{n-1}Q_{C_{k+1}}^{\,\,\,\,C_k}(z_{k+1,1},\ldots,z_{k+1,k+1};
\,z_{k,1},\ldots,z_{k,k}),
\ee
where $z_i:=z_{n,i}$ and the kernels $Q_{C_{k+1}}^{\,\,\,\,C_k}$
of the integral operators are given by the  convolutions
of the kernels  $Q_{C_{k+1}}^{\,\,\,\,D_{k+1}}$
and $Q_{D_{k+1}}^{\,\,\,\,C_k}$
 \bqa\label{QC}
Q_{C_{k+1}}^{\,\,\,\,C_k}(z_{k+1,i};\,z_{k,j})=\int\bigwedge_{i=1}^kdx_{k,i}
\, Q_{C_{k+1}}^{\,\,\,\,D_{k+1}}(z_{k+1,1},\ldots,z_{k+1,k+1};
x_{k+1,1},\ldots,x_{k+1,k+1})\times \\ \nonumber \times Q_{D_{k+1}}^{\,\,\,\,C_k}
(x_{k+1,1},\ldots,x_{k+1,k+1};z_{k,1},\ldots,z_{k,k}).\eqa
Thus the  wave function is given by the integral over a contour 
of the real dimension equal to a complex  dimension of the flag manifold
$X=G/B$, where $G=Sp(n,\mathbb{C})$ and $B$ is a Borel subgroup   
 \be
\sum_{k=1}^n (2k-1)=n^2=|R_+|.
\ee 

\subsection{$D_n$}

The eigenfunction for $D_n$ open Toda chain is given by
\be\label{WDN}
\Psi^{D_n}(x_1,\ldots,x_n)\,=\,
\int\bigwedge_{k=1}^{n-1}\bigwedge_{i=1}^kdx_{k,i}\,
\prod_{k=1}^{n-1}Q_{D_{k+1}}^{\,\,\,\,D_k}(x_{k+1,1},\ldots,x_{k+1,k+1};
\,x_{k,1},\ldots,x_{k,k}),\ee
where $x_i:=x_{n,i}$ and the kernels $Q_{D_{k+1}}^{\,\,\,\,D_k}$
of the integral operators are given
by the  convolutions  of the kernels $Q_{D_{k+1}}^{\,\,\,\,C_{k}}$
and $Q_{C_{k}}^{\,\,\,\,D_k}$
 \bqa
Q_{D_{k+1}}^{\,\,\,\,D_k}(x_{k+1,i};\,x_{k,j})=\int\bigwedge_{i=1}^kdz_{k,i}
\, Q_{D_{k+1}}^{\,\,\,\,C_k}(x_{k+1,1},\ldots,x_{k+1,k+1};
z_{k,1},\ldots,z_{k,k})\times \\ \nonumber \times Q_{C_{k}}^{\,\,\,\,D_k}
(z_{k,1},\ldots,z_{k,k};x_{k,1},\ldots,x_{k,k}).\eqa
Thus  the wave function is given by the integral over a contour 
of the real dimension equal to a complex  dimension of the flag manifold
$X=G/B$, where $G=SO(2n,\mathbb{C})$ and $B$ is a Borel subgroup   
 \be
\sum_{k=1}^{n-1}\, 2k=n(n-1)=|R_+|.
\ee

\section{ Givental diagrams and Toric degenerations}

In this section we describe combinatorial structure 
entering the integral representations
presented above.  This structure reflects 
flat toric degenerations  of the corresponding 
flag manifolds (see \cite{BCFKS} for $A_n$ and \cite{B}
for a  general approach to mirror symmetry via degeneration).
 The combinatorial structure 
of the integrand readily encoded into the (generalized) Givental
diagrams. Let us note that the diagrams for $B_n$, $C_n$ and $D_n$ 
can be obtained from those for $A_n$ by a factorization. This factorization
is closely related  with a particular realization \cite{DS} 
of  classical series of Lie algebras as fixed point subalgebras  of the
involutions acting on an algebra $\mathfrak{gl}_N$ for some $N$. 

\subsection{ $A_n$ Diagram }
Givental diagram \cite{Gi} for $A_n$ has the following form
\be\label{AnDiag}
\begin{CD}x_{n+1,n+1} @. @. @. @.\\ @AAb_{n,n}A @. @. @. @.\\
x_{n,n} @<<< x_{n+1,n} @. @. @.\\ @AAb_{n-1,n-1}A @AAb_{n,n-1}A @. @. @. @.\\
\vdots @. \vdots @. \ddots @.\\
@AAb_{2,2}A @AAb_{3,2}A @. @. @. @.\\ x_{2,2} @<a_{2,2}<< x_{3,2}
@<a_{3,2}<< \ldots @<a_{n,2}<< x_{n+1,2} @.\\ @AAb_{1,1}A
@AAb_{2,1}A @. @AAb_{n,1}A @.\\ x_{1,1} @<a_{1,1}<< x_{2,1}
@<a_{2,1}<< \ldots @<a_{n-1,1}<< x_{n,1} @<a_{n,1}<<
x_{n+1,1}\end{CD}\ee 
We assign   variables $x_{k,i}$ to the vertexes $(k,i)$ 
and  functions $e^{y-x}$  to the  arrows
 $(x\longrightarrow y)$ of the diagram
 (\ref{AnDiag}). The potential function $\mathcal{F}(x_{k,i})$ (see (\ref{pot}))
is given by the  sum of the functions assigned  to all arrows.

Note that the variables $\{x_{k,i}\}$ naturally parametrize 
an open part $U$ of the flag manifold $X=SL(n+1,\mathbb{C})/B$.   
The non-compact manifold $U$ has a natural action of the 
torus  and can be compactified to a (singular) toric variety. 
The set of the monomial relations defining this compactifiaction 
can be described as follows. Let us 
 introduce the new variables 
$$a_{k,i}=e^{x_{k,i}-x_{k+1,i}}
,\,\,\,\,\,b_{k,i}=e^{x_{k+1,i+1}-x_{k,i}},\,\,\,\,\,\,\quad 1\leq
k\leq n,\,\,1\leq i\leq k\,$$ assigned to the arrows of
the diagram (\ref{AnDiag}). Then the following 
defining relations hold  
\be\label{defrelAn}
 a_{k,i}\cdot b_{k,i}\,=\,
b_{k+1,i}\cdot a_{k+1,i+1},\qquad 1\leq k< n,\,\,1\leq i\leq k\\
a_{n,i}\cdot b_{n,i}=e^{x_{n,i+1}-x_{n,i}}\ee
The defining relations of the toric embedding  are given
by the monomial relations  for the  
variables associated with the paths on the diagram. They   
are given by a simple generalization of the relations
(\ref{defrelAn}) (see \cite{BCFKS} for details).

\subsection{$B_n$ Diagram}

Diagram for $B_n$ has the following form $(n=3)$ 
\be\begin{CD} @. @. @. @Vb_{31}VV @. @.\\
@. @. @>a_{31}>> z_{31} @>c_{31}>> x_{31} @. @.\\
@. @. @Vb_{21}VV @Vd_{21}VV @Vb_{32}VV @.\\
@. @>a_{21}>> z_{21} @>c_{21}>> x_{21} @>a_{32}>> z_{32} @>c_{32}>> x_{32} @.\\
@. @Vb_{11}VV @Vd_{11}VV @Vb_{22}VV @Vd_{22}VV @Vb_{33}VV @.\\
@>a_{11}>> z_{11} @>c_{11}>> x_{11} @>a_{22}>> z_{22} @>c_{22}>>
x_{22} @>a_{33}>> z_{33} @>c_{33}>> x_{33}\end{CD}\ee
Here we use the same rules for assigning 
variables to the arrows of the diagram as in $A_n$ case. In addition 
we assign functions $e^x$ to the arrows $(\longrightarrow
x)$.

Note that the diagram for $B_n$  can be obtained  by
factorization of the diagram (\ref{AnDiag}) for $A_{2n}$  by  the 
 following involution
\be \label{inv}
\iota\,:\quad X\longmapsto  w_0^{-1}X^T w_0,
\ee
where  $w_0$ is the longest element
of $A_{2n}$ Weyl group $\mathfrak{W}(A_{2n})$ isomorphic to a symmetric
group $\mathfrak{S}_{2n+1}$ and $X^T$ denotes the standard transposition.
Correspondingly the diagram for $B_n$ can be obtained from  $A_{2n}$ diagram  by
the quotient  
\be \label{w0invBn}
w_0\,:\qquad x_{k,i}\longleftrightarrow -x_{k,k+1-i}.
\ee
An analog of the monomial relations (\ref{defrelAn}) 
is as follows. 
Associate to the arrows of Givental diagram parameters
\be
\,a_{k,i}=e^{z_{k,i}-x_{k-1,i-1}},\,\,\,b_{k,i}=e^{z_{k,i}-x_{k,i-1}},\,\,
\,c_{k,i}=e^{z_{k,i}-x_{k,i}},\,\,\,d_{l,j}=e^{x_{l,j}-z_{l+1,j}}\,\,\,\\
\,1\leq k\leq n,\,\,1\leq i\leq k,
\quad 1\leq l\leq n-1,\,\, 1\leq j\leq l.\,\ee
Then the following relations hold: \bqa
a_{k,1}&=&b_{k,1},\qquad \qquad
\,\,\,\,\,\,\,\,\,\,\,\,\,\,\,\,\,\,\,\,1\leq k\leq n,\nonumber
\\
d_{k,i}\cdot a_{k+1,i+1}\,&=&\,c_{k+1,i}\cdot b_{k+1,i+1},\,\,\,\qquad 1\leq
k<n-1,\,\,1\leq i\leq k \\ b_{k,i}\cdot c_{k,i}\,&=&\,
a_{k+1,i}\cdot d_{k,i},\qquad \qquad \,1\leq k<n-1,\,\,1\leq i\leq k\nonumber\\
b_{n,i}\cdot c_{n,i}&=&e^{x_{n,i}-x_{n,i-1}}\nonumber\eqa

\subsection{ $C_n$ Diagram}
Diagram for $C_n$ has the following form $(n=3)$ 
\be \label{diagrCn}\begin{CD}
@. @.  z_{33} @. @. @. @.\\
@. @. @| @. @. @. @. @. @.\\
@. z_{22} @= x_{33} @<c_{33}<< z_{33} @.\\
@. @| @Ad_{22}AA @Ab_{32}AA @. @.\\
z_{11} @= x_{22} @<c_{22}<< z_{22} @<a_{32}<< x_{32} @<c_{32}<< z_{32} @.\\
 @| @Ad_{11}AA @Ab_{21}AA @Ad_{21}AA @Ab_{31}AA @. @.\\
x_{11} @<c_{11}<< z_{11} @<a_{21}<< x_{21} @<c_{21}<< z_{21}
@<a_{31}<< x_{31} @<c_{31}<< z_{31} @.\end{CD}\ee
where  one assigns functions 
$e^{-z-x}$ to  double arrows $(\begin{CD}x@= z\end{CD})$.

The Lie algebra $C_n$ can be realized as a fixed point subalgebra
of $A_{2n-1}$ using the involution
\be \label{inv}
\iota\,:\quad X\longmapsto  w_0^{-1}X^T w_0,
\ee
where  $w_0$ is the longest element of Weyl group
$\mathfrak{W}(A_{2n-1})=\mathfrak{S}_{2n}$ and $X^T$ denotes the standard transposition.
Correspondingly the diagram for $C_n$ can be obtained from  $A_{2n-1}$ diagram  by
the quotient  
\be \label{w0inv}
w_0\,:\qquad x_{k,i}\longleftrightarrow -x_{k,k+1-i}.
\ee
Note that diagram for $C_n$ can be also obtained 
by erasing the last row of vertexes
and arrows on the right slope from the diagram for
$D_{n+1}$ (see (\ref{D4Diag}) below).

An analog of the monomial relations (\ref{defrelAn})
is as follows. Let us introduce the  variables \be
a_{l,j}=e^{z_{l-1,j}-x_{l,j}},\qquad\qquad\qquad\,\,\,
\qquad 1<l\leq n,\,\,1\leq j\leq l,\,\, l\neq j,\\
b_{k,i}=e^{z_{k,i+1}-x_{k,i}},\,\,\,\qquad\qquad\qquad\qquad 1\leq k\leq n,\,\, 1\leq
i\leq k,\,\,k\neq i,\\
c_{k,i}=e^{x_{k,i}-z_{k,i}},\qquad \qquad\qquad \qquad
 1\leq k\leq n,\,\, 1\leq i\leq k,\\
d_{m,j}=e^{x_{m+1,j+1}-z_{m,j}},\qquad\qquad\qquad
\qquad 1\leq m<n,\,\,1\leq j\leq m\ee 
where one assign the variables $b_{11}$, $a_{22}$,
$b_{22}$, $a_{33}$, and $b_{33}$ to  the left slope of the diagram 
(\ref{diagrCn}). 
Then the following relations hold 
 \be c_{k,i}\cdot b_{k,i}\,=\,d_{k,i}\cdot
a_{k+1,i+1},\\ a_{k,i}\cdot d_{k-1,i}\,=\,b_{k,i}\cdot c_{k,i+1},\\
c_{n,i}\cdot b_{n,i}=e^{z_{n,i+1}-z_{n,i}},\qquad a_{n,n}\cdot
b_{n,n}=e^{-2z_{n,n}}.\ee

\subsection{ $D_n$ Diagram}

Diagram for $D_n$ has the following form $(n=4)$ 
\be\label{D4Diag}\begin{CD} \\
@. @.  z_{33} @= x_{44} @. @. @.\\
@. @. @| @Ad_{33}AA @. @. @. @. @. @.\\
@. z_{22} @= x_{33} @<c_{33}<< z_{33} @<a_{43}<< x_{43}\\
@. @| @Ad_{22}AA @Ab_{32}AA @Ad_{32}AA @.\\
z_{11} @= x_{22} @<c_{22}<< z_{22} @<a_{32}<< x_{32} @<c_{32}<<
z_{32} @<a_{42}<<
x_{42}\\
 @| @Ad_{11}AA @Ab_{21}AA @Ad_{21}AA @Ab_{31}AA @AAd_{31}A @.\\
x_{11} @<c_{11}<< z_{11} @<a_{21}<< x_{21} @<c_{21}<< z_{21}
@<a_{31}<< x_{31} @<c_{31}<< z_{31} @<a_{41}<< x_{41} \end{CD}\ee
Note that Lie algebra $D_n$ can be realized as a fixed point subalgebra
of $A_{2n-1}$ using the  involution
\be \label{inv}
\iota\,:\quad X\longmapsto  w_0^{-1}X^T w_0,
\ee
where  $w_0$ is the longest element of Weyl group
$\mathfrak{W}(A_{2n-1})=\mathfrak{S}_{2n}$ and $X^T$ denotes the standard transposition.
Correspondingly the  diagram for $D_n$ can be obtain
from that for $A_{2n-1}$ by the identification  of the
variables assigned to the vertexes of $A_{2n-1}$ diagram
\be \label{w0inv}
w_0\,:\qquad x_{k,i}\longleftrightarrow -x_{k,k+1-i}.
\ee

An analog of the monomial relations (\ref{defrelAn}) 
is as follows.  Let us introduce new variables 
associate to the arrows of the diagram 
 \be a_{l,j}=e^{z_{l-1,j}-x_{l,j}},\qquad\qquad\qquad\qquad\,\,\qquad\qquad\qquad
\qquad 1<l\leq n,\,\,1\leq j\leq l,\\
b_{k,i}=e^{z_{k,i+1}-x_{k,i}}\,\,\,c_{k,i}=e^{x_{k,i}-z_{k,i}},\,\,
\,d_{k,i}=e^{x_{k+1,i}-z_{k,i}},\qquad 1\leq k<n,\,\,1\leq i\leq k.\ee
The following defining relations hold
 \be c_{k,i}\cdot b_{k,i}\,=\,
d_{k,i}\cdot a_{k+1,i+1},\\ a_{k,i}\cdot d_{k-1,i}\,=\,b_{k,i}\cdot
c_{k,i+1},\\ a_{k,i}\cdot d_{k-1,i}=e^{x_{k,i+1}-x_{k,i}}.\ee 
Finally let us note that  it is easy to check that the potentials 
$\mathcal{F}(x,z)$ obtained from  $B_n$, $C_n$ and $D_n$ diagrams 
by summing the function assigned to the  arrows  
coincide with that entering  the integral representations in Section 4.   

\section{Elementary intertwiners for closed Toda chains}

In this section we generalize the construction of the elementary
intertwiners to the classical series of affine Lie algebras.
For the necessary facts in the theory of affine Lie algebras see
\cite{K},\,\cite{DS}.  Let us first recall the
construction of the $Q$-operator for $A_n^{(1)}$ closed Toda
chain~\cite{PG}.   The integral kernel of the  intertwining
$Q$-operator in this case reads \be
Q^{A_n^{(1)}}(x_1,\ldots,x_{n+1};\,y_1,\ldots,y_{n+1})=\exp\Big\{\,
\sum_{i=1}^{n+1}(e^{x_i-y_i}+g_{i+1}e^{y_{i+1}-x_i})\,\Big\},\quad
y_{n+2}=y_1,\ee 
The corresponding integral operator   intertwines the following 
Hamiltonians operators for $A_n^{(1)}$ closed Toda chains   \bqa
\CH^{A_n^{(1)}}(x_i)=-\frac{1}{2}\sum_{i=1}^{n+1}
\frac{\partial^2}{\partial
x_i^2}+g_1e^{x_1-x_{n+1}}+ \sum_{i=1}^ng_{i+1}e^{x_{i+1}-x_i},\\
\CH^{A_n^{(1)}}(y_i)=-\frac{1}{2}\sum_{i=1}^{n+1} \frac{\partial^2}{\partial
y_i^2}+g_1e^{y_1-y_{n+1}}+ \sum_{i=1}^ng_{i+1}e^{y_{i+1}-y_i}\,.\eqa

Below we provide  kernels of the  integral operators
intertwining Hamiltonians of the closed Toda chains corresponding 
to  other classical series of affine Lie algebras.

\subsection{ $A^{(2)}_{2n}\leftrightarrow BC_{n+1}^{(2)}$}

Simple roots  of  the twisted affine root system
$A^{(2)}_{2n}$ can be expressed  in terms of the standard basis $\{e_i\}$ as
follows \be  \alpha_1=e_1,\qquad \alpha_{i+1}=e_{i+1}-e_i,\qquad
1\leq i\leq n-1\qquad \alpha_{n+1}=-2e_n. \ee Corresponding
Dynkin diagram is given  by
$$\alpha_1
\Longleftarrow\alpha_2\longleftarrow\ldots\longleftarrow\alpha_{n-1}
\Longleftarrow\alpha_n.
$$
Simple roots of the twisted affine non-reduced  root system $BC^{(2)}_{n}$ are
 given by  \be\alpha_0=2e_1,\quad\alpha_1=e_1,\quad
\alpha_{i+1}=e_{i+1}-e_i,\,\,\, 1\leq i\leq n-1,\quad
\alpha_{n+1}=-e_n-e_{n-1}\ee and the corresponding Dynkin diagram is
as follows \bqa\begin{CD}\frac{\alpha_0}{\alpha_1}
\Longleftrightarrow\alpha_2@<<< \ldots@<<< \alpha_{n-1} @<<<
\alpha_n\\ @. @. @AAA @.\\ @. @. \alpha_{n+1} @.\end{CD}\eqa 
The integral operator with the following kernel \be \label{inttwA2even}
Q_{A^{(2)}_{2n}}^{BC^{(2)}_{n+1}}(x_i,\,z_i)=
\exp\Big\{\, g_1e^{z_1}+ \sum_{i=1}^n\Big(
e^{x_i-z_i}+g_{i+1}e^{z_{i+1}-x_i}\Big)+
g_{n+2}e^{-z_{n+1}-x_n}\,\Big\},\ee intertwines
Hamiltonian operators for $A_{2n}^{(2)}$ and $BC_{n+1}^{(2)}$ \bqa
\CH^{A_{2n}^{(2)}}(x_i)&=&-\frac{1}{2}\sum_{i=1}^{n}
\frac{\partial^2}{\partial x_i^2}+g_1e^{x_1}+
\sum_{i=1}^{n-1}g_{i+1}e^{x_{i+1}-x_i}+2g_{n+1}g_{n+2}e^{-2x_{n}},\\
\CH^{BC^{(2)}_{n+1}}(z_i)&=&-\frac{1}{2}\sum_{i=1}^{n+1} \frac{\partial^2}{\partial
z_i^2}+ \frac{g_1}{2}\Big(e^{z_1}+g_1e^{2z_1}\Big)+
\sum_{i=1}^{n}g_{i+1}e^{z_{i+1}-z_i}+g_{n+2}e^{-z_{n+1}-z_n}.\eqa
The integral kernel for the inverse transformation is  given by 
\be
Q^{\,\,\,\,A^{(2)}_{2n}}_{BC^{(2)}_{n+1}}(x_i,\,z_i)=Q_{A^{(2)}_{2n}}
^{\,\,\,\,BC^{(2)}_{n+1}}(z_i,\,x_i).\ee

\subsection{ $A^{(2)}_{2n-1}\leftrightarrow A^{(2)}_{2n-1}$}

Simple roots  of  the twisted affine root system
$A^{(2)}_{2n-1}$ are given by
\be \alpha_1=2e_1,\qquad \alpha_{i+1}=e_{i+1}-e_i,\qquad
1\leq i\leq n-1\qquad \alpha_{n+1}=-e_n-e_{n-1}, \ee and
corresponding Dynkin diagram is 
$$
\begin{CD}\alpha_1 \Longrightarrow \alpha_2 @>>> \ldots @>>> \alpha_{n-1}
@>>> \alpha_n\\ @. @. @VVV @.\\ @. @.\alpha_{n+1}@.\end{CD}
$$
The integral operator represented  by the following kernel
 \be\label{inttwA2odd}
Q_{A^{(2)}_{2n-1}}^{\,\,\,\,A^{(2)}_{2n-1}}(x_i,\,z_i)=\\
=\exp\Big\{\, g_1e^{x_1+z_1}+ \sum_{i=1}^{n-1}\Big(
e^{x_i-z_i}+g_{i+1}e^{z_{i+1}-x_i}\Big)+
e^{x_n-z_n}+g_{n+1}e^{-x_{n}-z_n}\,\Big\},\ee intertwines Hamiltonian
operators for $A^{(2)}_{2n-1}$ closed Toda chains with different
coupling constants \bqa
\CH^{A_{2n-1}^{(2)}}(x_i)&=&-\frac{1}{2}\sum_{i=1}^{n}
\frac{\partial^2}{\partial
x_i^2}+2g_1e^{2x_1}+\sum_{i=1}^{n-1}g_{i+1}e^{x_{i+1}-x_i}+
g_n g_{n+1}e^{-x_{n}-x_{n-1}},\\
\widetilde{\CH}^{A_{2n-1}^{(2)}}(z_i)&=&-\frac{1}{2}\sum_{i=1}^{n}
\frac{\partial^2}{\partial z_i^2}+ g_1g_2e^{z_1+z_2}+
\sum_{i=1}^{n-1}g_{i+1} e^{z_{i+1}-z_i}+2g_{n+1}e^{-2z_{n}}.\eqa

\subsection{$B^{(1)}_n\leftrightarrow BC^{(1)}_{n}$}

Simple roots of  the  affine root system
$B^{(1)}_{n}$ are given by 
 \be \alpha_1=e_1,\qquad \alpha_{i+1}=e_{i+1}-e_i,\qquad
1\leq i\leq n-1\qquad \alpha_{n+1}=-e_n-e_{n-1}, \ee and
corresponding Dynkin diagram is as follows
$$
\begin{CD}\alpha_1 \Longleftarrow \alpha_2 @<<< \ldots @<<< \alpha_{n-1}
@<<< \alpha_n\\ @. @. @AAA @.\\ @. @.\alpha_{n+1}@.\end{CD}
$$
Simple roots of  the affine non-reduced root system $BC^{(1)}_{n}$
are
\be \alpha_0=2e_1,\quad\alpha_1=e_1,\qquad
\alpha_{i+1}=e_{i+1}-e_i,\qquad 1\leq i\leq n-1\qquad
\alpha_{n+1}=-2e_n, \ee and corresponding Dynkin diagram is given by
$$
\frac{\alpha_0}{\alpha_1} \Longleftrightarrow \alpha_2 \longleftarrow \ldots
\longleftarrow \alpha_{n} \Longleftarrow\alpha_{n+1}$$
The integral operator represented by the following kernel
 \be\label{inttwBn}
Q_{B^{(1)}_n}^{\,\,\,\,\,BC^{(1)}_{n}}(x_i,\,z_i)=\\
=\exp\Big\{\, g_1e^{z_1}+ \sum_{i=1}^{n-1}\Big(
e^{x_i-z_i}+g_{i+1}e^{z_{i+1}-x_i}\Big)+
e^{x_n-z_n}+g_{n+1}e^{-x_{n}-z_n}\,\Big\},\ee intertwines
Hamiltonians of $B^{(1)}_n$ and $BC^{(1)}_{n}$ closed Toda chains
\bqa \CH^{B^{(1)}_n}(x_i)&=&-\frac{1}{2}\sum_{i=1}^{n}
\frac{\partial^2}{\partial
x_i^2}+g_1e^{x_1}+\sum_{i=1}^{n-1}g_{i+1}e^{x_{i+1}-x_i}+
g_ng_{n+1}e^{-x_{n}-x_{n-1}},\\
\CH^{BC^{(1)}_{n}}(z_i)&=&-\frac{1}{2}\sum_{i=1}^{n}
\frac{\partial^2}{\partial z_i^2}+ \frac{g_1}{2}\Big(e^{z_1}+g_1
e^{2z_1}\Big)+\sum_{i=1}^{n-1}g_{i+1}
e^{z_{i+1}-z_i}+2g_{n+1}e^{-2z_{n}}.\eqa
The integral kernel for the inverse transformation is  given by 
\be
Q^{\,\,\,\,B^{(1)}_{n}}_{BC^{(1)}_{n}}(x_i,\,z_i)=Q_{B^{(1)}_{n}}
^{\,\,\,\,BC^{(1)}_{n}}(z_i,\,x_i).\ee

\subsection{ $C^{(1)}\leftrightarrow D^{(1)}$}

Simple roots of  the affine root system
$C^{(1)}_{n}$ are
 \be \alpha_1=2e_1,\qquad \alpha_{i+1}=e_{i+1}-e_i,\qquad
1\leq i\leq n-1\qquad \alpha_{n+1}=-2e_n, \ee and corresponding
Dynkin diagram is given  by
$$
\alpha_1\Longrightarrow\alpha_2\longrightarrow\ldots
\longleftarrow\alpha_n\Longleftarrow\alpha_{n+1}
$$
Simple roots of  the affine root system
$D^{(1)}_{n}$ are
\be
\alpha_1=e_1+e_2,\qquad \alpha_{i+1}=e_{i+1}-e_i,\qquad 1\leq i\leq
n-1\qquad \alpha_{n+1}=-e_n-e_{n-1},
\ee
and corresponding Dynkin diagram is given  by
$$
\begin{CD}\alpha_1 @<<< \alpha_3 @<<< \ldots @<<< \alpha_{n-1}
@<<< \alpha_n\\ @. @AAA @. @AAA @.\\ @. \alpha_2 @. @. \alpha_{n+1}
@.\end{CD}
$$
The integral operator  with the following kernel \be\label{inttwCn}
Q_{C^{(1)}_{n}}^{\,\,\,\,D^{(1)}_{n+1}}(x_i,\,z_i)=
\exp\Big\{ g_1e^{x_1+z_1}+\sum_{i=1}^{n}\Big(
e^{x_i-z_i}+ g_{i+1}e^{z_{i+1}-x_i}\Big)+
g_{n+2}e^{-z_{n+1}-x_n}\Big\}\ee
 intertwines Hamiltonian operators for $C^{(1)}_n$ and $D^{(1)}_{n+1}$
 closed Toda chains
\bqa
\CH^{C^{(1)}_{n}}(x_i)&=&-\frac{1}{2}\sum_{i=1}^{n}
\frac{\partial^2}{\partial x_i^2}+2g_1e^{2x_1}+
\sum_{i=1}^{n-1}g_{i+1}e^{x_{i+1}-x_i}+
2g_{n+1}g_{n+2}e^{-2x_{n}},\\
 \CH^{D^{(1)}_{n+1}}(z_i)&=&-\frac{1}{2}\sum_{i=1}^{n+1}
\frac{\partial^2}{\partial z_i^2}+g_1g_2e^{z_1+z_2}+
\sum_{i=1}^{n}g_{i+1}e^{z_{i+1}-z_i}+g_{n+2}e^{-z_{n+1}-z_{n}}.\eqa
The integral operator with the kernel
 \be\label{inttwDn}
Q_{D^{(1)}_{n}}^{\,\,\,\,C^{(1)}_{n-1}}(x_i,\,z_i)=\exp\Big\{\,
g_1e^{x_1+z_1}+ \sum_{i=1}^{n-1}\Big(e^{z_i-x_i}+
g_{i+1}e^{x_{i+1}-z_i}\Big)+ g_{n+1}e^{-x_{n}-z_{n-1}}\,\Big\}\ee
 intertwines Hamiltonian operators for $D^{(1)}_n$ and $C^{(1)}_{n-1}$
 closed Toda chains
 \bqa
\CH^{D^{(1)}_{n}}(x_i)&=&-\frac{1}{2}\sum_{i=1}^{n} \frac{\partial^2}{\partial
x_i^2}+g_1g_2e^{x_1+x_2}+
\sum_{i=1}^{n-1}g_{i+1}e^{x_{i+1}-z_i}+g_{n+1}e^{-x_{n}-x_{n-1}},\\
\CH^{C^{(1)}_{n-1}}(z_i)&=&-\frac{1}{2}\sum_{i=1}^{n-1} \frac{\partial^2}{\partial
z_i^2}+2g_1e^{2z_1}+ \sum_{i=1}^{n-2}g_{i+1}e^{z_{i+1}-z_i}+
2g_{n}g_{n+1}e^{-2z_{n-1}}.\eqa
The integral kernel for the inverse transformation is given by 
\be
Q^{\,\,\,\,D^{(1)}_{n}}_{C^{(1)}_{n}}(x_i,\,z_i)=Q_{D^{(1)}_{n}}
^{\,\,\,\,C^{(1)}_{n}}(z_i,\,x_i).\ee

\section{Baxter $Q$-operators}

Now we apply the results presented in the previous Sections to  the
construction of the  Baxter integral $Q$-operators for all classical series
 of affine Lie algebras. Let us note that the elementary intertwining
operators and recursive operators for finite Lie algebras can be obtained
from the elementary intertwining operators and Baxter operators
for affine Lie algebras by taking appropriate limits $g_i\rightarrow0$ in
(\ref{inttwA2even})-(\ref{inttwDn}). This generalizes  known
relation between Baxter operator for $A_n^{(1)}$ and recursive
operators for $A_n$.

The integral kernels for $Q$-operators have the following form
\bqa
Q^{A^{(2)}_{2n}}(x_1,\ldots,x_n;\,y_1,\ldots,y_n)=\int
\bigwedge_{i=1}^{n+1}dz_i
Q_{A^{(2)}_{2n}}^{\,\,\,\,BC^{(2)}_{n+1}}(x_1,\ldots,x_n;\,z_1,\ldots,z_{n+1})\times
\\ \nonumber
\times
Q_{BC^{(2)}_{n+1}}^{\,\,\,\,A^{(2)}_{2n}}(z_1,\ldots,z_{n+1};\,y_1,\ldots,y_n),\eqa
\bqa
Q^{A^{(2)}_{2n-1}}(x_1,\ldots,x_n;\,y_1,\ldots,y_n)=\int
\bigwedge_{i=1}^{n}dz_i\,
Q_{A^{(2)}_{2n-1}}^{\,\,\,\,A^{(2)}_{2n-1}}(x_1,\ldots,x_n;\,z_1,\ldots,z_{n})\times
\\ \nonumber
\times Q_{A^{(2)}_{2n-1}}^{\,\,\,\,
A^{(2)}_{2n-1}}(y_1,\ldots,y_{n};\,z_1,\ldots,z_n),\eqa
\bqa
Q^{B^{(1)}_{n}}(x_1,\ldots,x_n;\,y_1,\ldots,y_n)=\int
\bigwedge_{i=1}^{n}dz_i\,
Q_{B^{(1)}_{n}}^{\,\,\,\,BC^{(1)}_{n}}(x_1,\ldots,x_n;\,z_1,\ldots,z_{n})
\times \\ \nonumber
Q_{BC^{(1)}_{n}}^{\,\,\,\,B^{(1)}_{n}}(z_1,\ldots,z_{n};\,y_1,\ldots,y_n),\eqa
\bqa
Q^{C^{(1)}_{n}}(x_1,\ldots,x_n;\,y_1,\ldots,y_n)=\int
\bigwedge_{i=1}^{n+1}dz_i\,
Q_{C^{(1)}_{n}}^{\,\,\,\,D^{(1)}_{n+1}}(x_1,\ldots,x_n;\,z_1,\ldots,z_{n+1})
\times \\ \nonumber
Q_{D^{(1)}_{n+1}}^{\,\,\,\,C^{(1)}_{n}}(z_1,\ldots,z_{n+1};\,y_1,\ldots,y_n),\eqa
\bqa
Q^{D^{(1)}_{n}}(x_1,\ldots,x_n;\,y_1,\ldots,y_n)=\int
\bigwedge_{i=1}^{n-1}dz_i\,
Q_{D^{(1)}_{n}}^{\,\,\,\,C^{(1)}_{n-1}}(x_1,\ldots,x_n;\,z_1,\ldots,z_{n-1})
\times \\ \nonumber
Q_{C^{(1)}_{n-1}}^{\,\,\,\,D^{(1)}_{n}}(z_1,\ldots,z_{n-1};\,y_1,\ldots,y_n).\eqa

\section{Baxter operators for $B_{\infty}$,  $C_{\infty}$ and
$D_{\infty}$}

Similar approach can be applied to construct  Baxter $Q$  
operators for infinite root systems $B_{\infty}$, $C_{\infty}$, $BC_{\infty}$ and
$D_{\infty}$.  

Simple roots and Dynkin diagrams for infinite Lie algebras 
$A_{\infty}$, $B_{\infty}$, $C_{\infty}$, $BC_{\infty}$  and
$D_{\infty}$ are as follows 
\be A_{\infty}:\,\,\,\,\,\,\,\,\,\,\,\,\,\qquad\qquad\,\,
\qquad  \alpha_{i+1}=e_{i+1}-e_i, \qquad \quad
i\in\mathbb{Z}, \ee \be\begin{CD}
\ldots @>>> \alpha_{-1}\,@>>>\,\alpha_0 @>>> \alpha_1 @>>> \ldots
\end{CD}\ee

\be B_{\infty}:\,\,\,\,\,\,\,\,
\alpha_1=e_1,\qquad \alpha_{i+1}=e_{i+1}-e_i, \qquad\qquad
i\in\mathbb{Z}_{>0}, \ee \be\begin{CD}
\alpha_1\,\Longleftarrow\,\alpha_2 @<<< \alpha_3 @<<< \ldots
\end{CD}\ee

\be C_{\infty}:\,\,\,\,\,\,\, \alpha_1=2e_1,\qquad
\alpha_{i+1}=e_{i+1}-e_i, \qquad \qquad i\in\mathbb{Z}_{>0}, \ee
\be\begin{CD} \alpha_1\,\Longrightarrow\,\alpha_2 @>>> \alpha_3 @>>>
\ldots \end{CD}\ee

\be
D_\infty:\,\,\,\,\,\,\
\alpha_1=e_1+e_2,\qquad \alpha_{i+1}=e_{i+1}-e_i, \qquad\quad i\in\mathbb{Z}_{>0},
\ee
\be
\begin{CD}\alpha_1 @>>> \alpha_3 @>>> \alpha_4
@>>> \ldots\\ @. @AAA @. @.\\ @. \alpha_2 @. @.
\end{CD}
\ee

\be BC_{\infty}:\,\,\,\,\,\, \alpha_0=2e_1,\qquad
\alpha_1=e_1,\qquad \alpha_{i+1}=e_{i+1}-e_i, \qquad
i\in\mathbb{Z}_{>0}, \ee
\be\begin{CD}\frac{\alpha_0}{\alpha_1}\Longleftrightarrow \alpha_2
@<<< \alpha_3 @<<<\ldots\end{CD}\ee 
Baxter $Q$-operator for $A_\infty$ infinte Toda chain is known 
(see \cite {T} for the classical limit).
It is given by an  integral operator with the kernel 
 \be
Q(x_i,\,y_i)=\exp\Big\{\,\sum_{i\in\mathbb{Z}}\left(
e^{x_i-y_i}+g_ie^{y_{i+1}-x_i}\right)\,\Big\},\ee and intertwines
$A_{\infty}$ Toda chain Hamiltonians
 \bqa
\CH^{A_{\infty}}(x_i)&=&-\frac{1}{2}\sum_{i\in\mathbb{Z}}\frac{\partial^2}{\partial
x_i^2}+\sum_{i\in \mathbb{Z}} g_i e^{x_{i+1}-x_i},\\
\CH^{A_{\infty}}(y_i)&=&-\frac{1}{2}\sum_{i\in\mathbb{Z}}\frac{\partial^2}{\partial
y_i^2}+\sum_{i\in \mathbb{Z}}g_ie^{y_{i+1}-y_i}.\eqa
Its generalization to  other classical series relies on  
 the construction of the integral  operators
intertwining different classical series. Thus we have the following
set of integral operators. 

The integral operator with the kernel \be
Q_{B_\infty}^{\,\,\,\,\,BC_\infty}(x_i,\,z_i)=\exp\Big\{\,
g_1e^{z_1}+ \sum_{i>0}
\left(e^{x_i-z_i}+g_{i+1}e^{z_{i+1}-x_i}\right)\,\Big\},\ee intertwines $B_\infty$
and $BC_\infty$ Toda chain Hamiltonian operators 
 \bqa
\CH^{B_{\infty}}(x_i)&=&-\frac{1}{2}\sum_{i=1}^\infty\frac{\partial^2}{\partial x_i^2}
+g_1e^{x_1}+\sum_{i=1}^\infty g_{i+1}e^{x_{i+1}-x_i},\\
\CH^{BC_{\infty}}(z_i)&=&-\frac{1}{2}\sum_{i=1}^\infty\frac{\partial^2}{\partial z_i^2}
+\frac{g_1}{2}\Big(e^{z_1}+g_1e^{2z_1}\Big)+\sum_{i=1}^\infty
g_{i+1}e^{z_{i+1}-z_i}.\eqa
$Q$-operator for $B_\infty$ is then obtained
by composition of the intertwiner operators. 
For the integral kernel of the $Q$-operator for $B_{\infty}$ Toda
chain we have 
$$Q^{B_\infty}(x_i;y_i)=\int
\bigwedge_{i=1}^\infty dz_i\, Q_{B_\infty}^{\,\,\,\,BC_\infty}(x_i,\,z_i)
\cdot Q_{BC_\infty}^{\,\,\,\,\,B_\infty}(z_i,\,y_i).$$
Similarly the integral operator with the kernel
 \be
Q_{C_\infty}^{\,\,\,\,D_\infty}(x_i,\,z_i)=
\exp\Big\{\,g_1e^{x_1+z_1}+
\sum_{i>0} \left(e^{x_i-z_i}+g_{i+1}e^{z_{i+1}-x_i}\right)\,\Big\},\ee intertwines
$C_\infty$ and $D_\infty$  Toda chain Hamiltonian operators
 \bqa
\CH^{C_{\infty}}(x_i)&=&-\frac{1}{2}\sum_{i=1}^\infty\frac{\partial^2}{\partial
x_i^2} +2g_1e^{2x_1}+\sum_{i=1}^\infty g_{i+1}e^{x_{i+1}-x_i},\\
\CH^{D_{\infty}}(z_i)&=&-\frac{1}{2}\sum_{i=1}^\infty\frac{\partial^2}{\partial z_i^2}
+g_1g_2e^{z_1+z_2}+\sum_{i=1}^\infty g_{i+1}e^{z_{i+1}-z_i}.\eqa
Thus integral $Q$-operator for  $C_\infty$ has as the  kernel
$$Q^{C_\infty}(x_i;y_i)=\int \bigwedge_{i=1}^\infty dz_i\,
Q_{C_\infty}^{\,\,\,\,\,D_\infty}(x_i,\,z_i) \cdot
Q_{D_{\infty}}^{\,\,\,\,C_{\infty}}(z_i,\,y_i).$$ 
The integral operator with the kernel \be
Q_{D_\infty}^{\,\,\,\,C_\infty}(z_i,\,x_i)=\exp\Big\{\,
g_1e^{x_1+z_1}+ \sum_{i>0}\left(
e^{z_i-x_i}+g_{i+1}e^{x_{i+1}-z_i}\right)\,\Big\},\ee intertwines
$D_{\infty}$ and $C_{\infty}$ Toda chain Hamiltonian operators \bqa
\CH^{D_{\infty}}(x_i)&=&-\frac{1}{2}\sum_{i=1}^\infty\frac{\partial^2}{\partial
x_i^2} +g_1g_2e^{x_1+x_2}+\sum_{i=1}^\infty g_{i+1}e^{x_{i+1}-x_i},\\
\CH^{C_{\infty}}(z_i)&=&-\frac{1}{2}\sum_{i=1}^\infty\frac{\partial^2}{\partial z_i^2}
+2g_1e^{2z_1}+\sum_{i=1}^\infty g_{i+1}e^{z_{i+1}-z_i}.\eqa Therefore
 Baxter integral  $Q$-operator for  $D_\infty$ has the following
 kernel
$$
Q^{D_\infty}(x_i;y_i)=\int \bigwedge_{i=1}^\infty dz_i\cdot
Q_{D_\infty}^{\,\,\,\,C_\infty}(x_i,\,z_i) \cdot
Q_{C_\infty}^{\,\,\,\,D_\infty}(z_i,\,y_i).$$

\section{Conclusions}

In this note  we provide explicit  expressions for recursive operators
 and Baxter $Q$-operators for classical series of Lie algebras.
This allows us to generalize Givental representation for the
eigenfunctions of the  open Toda chain
quadratic Hamiltonians to  all classical series. 
In this note we consider only the case of zero eigenvalues
leaving a rather   straightforward
 generalization to  non-zero eigenvalues to another occasion\cite{GLO}.
The proof of the eigenfunction property 
for the full set of  Toda chain Hamiltonians 
follows the same strategy as in \cite{Gi},\,\cite{GKLO}
 and will be published separately.
The results presented in this note provide 
 a generalization to other  classical series 
of only a part of \cite{GKLO}. The other part connected with the
 interpretation of  elementary intertwining operators in
 representation theory framework  will discussed elsewhere.

Let us note that $Q$-operator was introduced by Baxter
as a key tool to  solve quantum integrable systems (see \cite{B}).
Therefore one can expect that the  Givental representation
and its generalizations should play an important 
role in the  theory of quantum integrable systems
solved by the quantum inverse scattering method 
(see e.g.  \cite{Fa}).

Finally let us mention two possible applications of the 
obtained results. The generalization  of the Givental integral 
representation and corresponding  diagrams to other classical 
series  allows to describe  flat torifications 
of the corresponding flag manifolds. This provides interesting
applications to the explicit description of 
the  mirror symmetry for flag manifolds. In particular 
one expects the mirror/Langlands duality between Gromov-Witten 
invariants for flag manifolds and the corresponding period integrals
 for $B_n$ and $C_n$ series. 
Another interesting  application connected with the 
theory of automorphic forms and corresponding $L$-functions.  
Note that  open Toda chain wave functions are given (see
e.g. \cite{STS}) by the Whittaker functions \cite{Ha} for the
corresponding Lie groups. 
Taking into account the simplicity and the unified form of
the arising constructions of Whittaker functions 
for all  classical groups one can  expect important
computational advantages in using the proposed 
integral representations.

\end{document}